\documentclass{article}
\usepackage{spconf,amsmath,graphicx} % Required for inserting images

\usepackage{cite}
\usepackage{amsmath,amssymb,amsfonts}
\usepackage{algorithm2e}
\usepackage{graphicx}
\usepackage{textcomp}
\usepackage{xcolor}
\usepackage{subcaption}
\usepackage{url}
\usepackage{bm}
\usepackage{alphabeta}

\newcommand{\bA}{{\mathbf{A}}}
\newcommand{\bI}{{\mathbf{I}}}
\newcommand{\bx}{{\bm{x}}}
\newcommand{\condA}{{\mathrm{cond}(\bA)}}

\begin{document}
\name{Chai Wah Wu, Mark S. Squillante, Vasileios Kalantzis, and Lior Horesh}
\address{IBM Research, Thomas J. Watson Research Center, USA}

\ninept

\title{Stable iterative refinement algorithms for solving linear systems}

\maketitle
\begin{abstract}
    Iterative refinement (IR) is a popular scheme for solving a linear system of equations based on gradually improving the accuracy of an initial approximation. Originally developed to improve upon the accuracy of Gaussian elimination, interest in IR has been revived because of its suitability for execution on fast low-precision hardware such as analog devices and graphics processing units. IR generally converges when the error associated with the solution method is small, but is known to diverge when this error is large. We propose and analyze a novel enhancement to the IR algorithm by adding a line search optimization step that guarantees the algorithm will not diverge. Numerical experiments verify our theoretical results and illustrate the effectiveness of our proposed scheme.
\end{abstract}


\section{Introduction}

Many problems in signal processing and engineering are concerned with the task of approximately decomposing a vector representation of data through 
a linear combination of fixed basis signals \cite{markovsky2007overview,cadzow1990signal,selesnick2013least}. Such problems can be typically solved via Linear Least Squares (LLS), a workhorse algorithm in signal processing employed in commonly used applications such as causal Wiener filters, system identification, linear prediction, and filter 
design \cite{markovsky2007overview,wang2011least,aastrom1971system}. 
At the core of LLS lies the minimization of $\|\mathbf{H}\bm{x}-\bm{y}\|^2$ where $\mathbf{H}$ and $\bm{y}$ are given. For example, in massive Multiple Input Multiple Output (MIMO) applications, the matrix 
$\mathbf{H}$ denotes an $\Mu\times \Nu$ channel matrix encoding the uplink relationship between $\Mu$ receiving antennas at a base station and $\Nu$ single-user antennas, and the vector $\bm{y}$ denotes 
the received vector subject to additive white Gaussian error \cite{lu2014overview,marzetta2015massive,larsson2014massive}.

When $\mathbf{H}$ is a full-rank matrix, the solution of the LLS problem can be computed 
via the \emph{normal equations} $\bm{x} = (\mathbf{H}^T\mathbf{H})^{-1}\mathbf{H}^T\bm{y}$.
The computation of $\bm{x}$ is  achieved by first computing the product $\bm{b} = 
\mathbf{H}^T\bm{y}$ and then solving 
the symmetric and positive-definite (SPD) linear system $\mathbf{A}\bm{x}=\bm{b}$ where $\mathbf{A} = \mathbf{H}^T\mathbf{H}$. The cost to form the matrix $\mathbf{A}$ is asymptotically equal to $O(\Mu \Nu^2)$ and the corresponding linear 
system can be solved by computing the Cholesky decomposition $\mathbf{A}=\mathbf{L}\mathbf{L}^T$ and performing 
a triangular substitution with the matrices $\mathbf{L}$ and $\mathbf{L}^T$.
The cost of the Cholesky decomposition 
is equal to $O(\Nu^3)$, which can be impractical for large $\Nu$. 

An alternative option is to 
solve $\mathbf{A}\bm{x}=\bm{b}$ approximately via an iterative method. In contrast to 
direct methods, iterative methods skip the matrix decomposition step and instead 
refine an initial approximation of $\bm{x}$ so that, hopefully, a good estimate can be produced after a small number of iterations \cite{saad2003iterative}. Normal equations is just one application where the solution of linear systems 
becomes necessary in signal processing, e.g., see  \cite{zhao2000hierarchical,chang2000analysis,tu2020efficient,wang2022randomized} for a non-exhaustive list.
When $\mathbf{H}$ is a full-rank matrix, the normal equations matrix 
$\mathbf{A}$ is symmetric and positive-definite, and the linear system 
$\mathbf{A}\bm{x}=\bm{b}$ can be solved via the Conjugate Gradient (CG) method \cite{saad2003iterative}. Nonetheless, the condition number of the matrix 
$\mathbf{A}$ is approximately equal to the square of that of $\mathbf{H}$, which can 
lead to an ill-conditioned system and slow convergence in CG as well as Krylov subspace methods in general \cite{liesen2013krylov}. 

To remedy the solution of potentially ill-conditioned linear systems, Wilkinson~\cite{Wilkinson} suggested the Iterative Refinement (IR) 
algorithm, an iterative method which during the $m$-th iteration corrects 
the approximation $\bm{x}_m$ of $\bm{x}$ to the enhanced 
approximation $\bm{x}_{m+1}$. IR is based on an inner-outer 
iteration, where the outer scheme maintains control of the iterative 
procedure by computing the residual of the approximation whereas the 
inner iteration performs a linear system solution with the matrix $\mathbf{A}$ that might not be highly accurate. Except for use in the solution of 
ill-conditioned linear systems, inner-outer iterative schemes 
such as IR have been applied in matrix preconditioning \cite{golub1999inexact} and PageRank \cite{gleich2010inner}, 
as well as high-performance computing applications (HPC) \cite{haidar2017investigating}.

\vspace{0.05in}

\noindent\textbf{Contribution:} While IR can lead to fast convergence towards a highly accurate approximation of $\bm{x}$, IR might diverge when $\mathbf{A}$ is ill-conditioned or the inner iteration introduces more error than the outer iteration can sustain. The main contribution of this paper is 
an enhanced IR scheme that does not result in divergence regardless 
of the magnitude of the error stemming from the approximate linear system solutions or the condition number of the matrix $\mathbf{A}$. More specifically, we present a new theoretical approach that introduces a line search step along the newly computed direction. Our analysis indicates that multiplying the search direction with the scalar which minimizes the line search objective is guaranteed to not magnify the error. We illustrate that the proposed IR framework can converge faster and converges even when standard IR fails to converge. Moreover, we show that our enhanced IR scheme can converge even when the inner iteration relies on fast but noisy analog hardware.

\vspace{0.05in}

\noindent\textbf{Notation:} Greek letters denote scalars, lowercase bold 
letters, e.g., $\bm{x}$, denote vectors, and uppercase bold letters, e.g., 
$\mathbf{A}$, represent matrices. Moreover, the operator $\condA$ denotes 
the condition number of the matrix $\bA$. The
$(i,j)$-th entry of matrix $\mathbf{A}$ is denoted by
$[\mathbf{A}]_{ij}$. 
Finally, $\bI$ denotes the identity matrix of the appropriate dimension. 
Unless specified otherwise, the norm $\|.\|$ represents the Euclidean 
norm of a vector and the induced 2-norm of a matrix.

\section{Iterative refinement}\label{sec:ir}
The original form of the IR scheme, first proposed by Wilkinson~\cite{Wilkinson} and further analyzed by Moler~\cite{Moler}, uses a low precision but fast method, referred to as a {\em basic} method, together with full precision arithmetic to iteratively reduce the errors in solving a (dense) linear system of equations $\mathbf{A}\bm{x}=\bm{b}$,
especially when the matrix $\mathbf{A}$ is ill-conditioned.
The main idea is that an inaccurate solver is used in the basic method to do the heavy lifting and solve the residual equation $\mathbf{A}\bm{d} = \bm{r} \triangleq \bm{b}-\mathbf{A}\bm{x}$, requiring $O(\Nu^3)$ operations, and to then update the estimate $\bm{x}$ and compute the next residual requiring $O(\Nu^2)$ operations,
with the residual computed using higher precision arithmetic.
It was shown by Wilkinson and Moler that if
$\condA$
is not too large and the residual is computed at double the working precision, then IR will converge to the true solution within working precision. The standard IR algorithm is presented as Algorithm \ref{alg:ir}. 

\RestyleAlgo{ruled}
\begin{algorithm}
\caption{Iterative Refinement}
\label{alg:ir}
\KwData{$\mathbf{A}$, $\bm{b}$, $\bm{x}_0$}
\KwResult{$\bm{x}$ such that $\mathbf{A}\bm{x}=\bm{b}$}
\While{Stopping criteria is not satisfied}{
compute residual $\bm{r}_m \gets \bm{b}-\mathbf{A}\bm{x}_m$\;
solve for $\bm{z}$ in $\mathbf{A}\bm{z} = \bm{r}_m$ using the basic method resulting in a solution $\bm{d}_m$\;
update $\bm{x}_{m+1} \gets \bm{x}_m + \bm{d}_m$\;
}
\end{algorithm}

In recent years, the ever increasing peak performance of fast low-precision hardware, e.g., Graphics Processing Units (GPUs) or analog crossbar arrays, has 
motivated the use of such hardware to perform the reduced precision step at 
each IR iteration \cite{haidar2018harnessing,haidar2020mixed,haidar2018design,kalantzis2021solving,kalantzis2023fgmres}. 
The use of low-precision, possibly non-digital, hardware leads to computing substrates such as stochastic rounding or dither rounding, which allow power, 
speed or size efficiency, but incur runtime numerical errors that can be either deterministic or stochastic, and of much higher magnitude compared to conventional 
digital hardware executing in single precision \cite{stochastic_computing_survey,dither_computing}. 
Thus, leveraging the benefits of such hardware requires stable IR variants that introduce safeguard steps in order to avoid divergence 
during the iterative procedure.

\section{Stable iterative refinement}\label{sec3}
It has been observed that
IR
can diverge when the error in the basic method to solve $\mathbf{A}\bm{d}_m = \bm{r}_{m}$ is large \cite{Wilkinson,Moler}. 
Specifically, letting $\bm{x} = \mathbf{A}^{-1}\bm{b}$ be the true solution of the problem, IR does not guarantee that $\|\bm{x}_{m+1}-\bm{x}\| \leq \|\bm{x}_{m}-\bm{x}\|$;
e.g., 
this is true when the Euclidean norm is used, $\mathbf{A}$ is nonsingular, and $\bm{r}_{m}$ is in the opposite direction as $\mathbf{A}^{-T}\bm{d}_{m}$. In particular, if $\bm{r}_{m}^T\mathbf{A}^{-T}\bm{d}_{m} \leq 0$, then 
\begin{align*}
    \bm{r}_{m}^T\mathbf{A}^{-T}\bm{d}_{m} &= (\bm{b}-\mathbf{A}\bm{x}_{m})^T\mathbf{A}^{-T}\bm{d}_{m}
    = (\bm{x}-\bm{x}_{m})^T\bm{d}_{m} 
    = 0,
    \\ \mbox{i.e.,} \qquad\qquad\quad & \\
    \|\bm{x}_{m+1}-\bm{x}\|^2 &= \|\bm{x}_{m}-\bm{x}\|^2 + \|\bm{d}_{m}\|^2 - 2\bm{r}_{m}^T\mathbf{A}^{-T}\bm{d}_{m}. 
\end{align*}
Similarly, when $\bm{r}_{m}$ is in the opposite direction as $\mathbf{A}\bm{d}_m$, then $\|\bm{r}_{m+1}\|$ can be larger than $\|\bm{r}_{m}\|$. Specifically,  
\begin{equation*}
   \|\bm{r}_{m+1}\|^2 = \|\bm{r}_{m} - \mathbf{A}\bm{d}_m\|^2 = \|\bm{r}_{m}\|^2 + \|\mathbf{A}\bm{d}_m\|^2 - 2\bm{r}_{m}^T\mathbf{A}\bm{d}_m ,
\end{equation*}
and thus $\|\bm{r}_{m+1}\|^2 > \|\bm{r}_{m}\|^2$ when $\bm{r}_{m}^T\mathbf{A}\bm{d}_m < 0$.

To correct this problem, we replace the update equation in IR with a line search along the direction $\bm{d}_{m}$ to minimize $\|\bm{b} - \mathbf{A}(\bm{x}_{m} + \alpha \bm{d}_{m})\|$, resulting in a stable IR algorithm provided in Algorithm \ref{alg:sir}. 
The additional computational complexity of our stable IR algorithm over
the standard IR algorithm is dominated by the addition of only two vector inner products (to compute $\alpha_m$) which represents a small fraction of the overall computation of both methods.

\RestyleAlgo{ruled}
\begin{algorithm}
\caption{Stable Iterative Refinement}
\label{alg:sir}
\KwData{$\mathbf{A}$, $\bm{b}$, $\bm{x}_0$}
\KwResult{$\bm{x}$ such that $\mathbf{A}\bm{x}=\bm{b}$}
compute residual $\bm{r}_0 \gets \bm{b}-\mathbf{A}\bm{x}_0$\;
\While{Stopping criteria is not satisfied}{
solve for $\bm{z}$ in $\mathbf{A}\bm{z} = \bm{r}_{m}$ using the basic method resulting in a solution $\bm{d}_{m}$\;
compute $\bm{w}_m=\mathbf{A}\bm{d}_m$\;
compute $\alpha_m=\bm{r}_{m}^T\bm{w}_m/\|\bm{w}_m\|^2$\;
update $\bm{x}_{m+1} \gets \bm{x}_{m} + \alpha_m \bm{d}_{m}$\;
update $\bm{r}_{m+1} \gets \bm{r}_{m}-\alpha_m \bm{w}_m$\;
}
\end{algorithm}

\subsection{Nondivergence of stable IR}

We show that Algorithm \ref{alg:sir} does not diverge regardless of how inaccurate or random the basic method is.
In particular, we have 
\begin{align*}
    \alpha_m &= \mathrm{argmin}_\alpha \|\bm{b} - \mathbf{A}(\bm{x}_{m} + \alpha \bm{d}_{m})\|
    \\ 
    &
    = \mathrm{argmin}_\alpha \|\bm{r}_{m}- \alpha \mathbf{A}\bm{d}_m\|, 
\end{align*}
with $\bm{r}_{m+1}=\bm{r}_{m}- \alpha_m \mathbf{A}\bm{d}_m$ and therefore $\|\bm{r}_{m+1}\| \leq \|\bm{r}_{m}\|$, which means that
 $\|\bm{r}_{m}\|$ converges. If $\mathbf{A}$ is nonsingular, then this implies that $\|\bm{x}_{m} - \bm{x}\|$ converges since $\bm{r}_{m} = \mathbf{A}(\bm{x}-\bm{x}_{m})$.
 Our desired goal is that $\|\bm{r}_{m}\| \rightarrow 0$ and thus $\bm{x}_{m}\rightarrow \bm{x}$.

 Furthermore, starting from the same initial $\bm{r}_{m}$, our stable IR algorithm generates $\bm{r}_{m+1}$ that has smaller or equal norm than the $\bm{r}_{m+1}$ generated by IR. This does not mean that, starting from the same initial $\bm{x}_{m}$, the modified IR algorithm generates $\bm{x}_{m+1}$ that is closer to $\bm{x}$ than the $\bm{x}_{m+1}$ generated by IR.
 However, we show in the next section that the convergence criteria for IR with respect to $\bm{x}_{m}$ is similar to that of our modified IR.
In other words, if
$\condA$
and the error of the basic method are small enough, then our modified IR will also converge. 
Note that if we set $\alpha_m = \bm{d}_{m}^T(\bm{x}-\bm{x}_{m})\|\bm{d}_{m}\|^{-2}$, then we can guarantee that $\|\bm{x}_{m+1}-\bm{x}\| \leq \|\bm{x}_{m}-\bm{x}\|$.
Unfortunately, this computation is not realistic since it involves the unknown  $\bx$.

\subsection{Convergence of stable IR}
The basic method typically solves $\mathbf{A}\bm{z}=\bm{r}_{m}$ inaccurately, i.e., $\mathbf{A}\bm{d}_m\neq \bm{r}_{m}$. Supposing we can write $\bm{r}_{m}-\mathbf{A}\bm{d}_m$ as $\mathbf{A}\mathbf{F}_{m}\bm{d}_{m}$, it follows that 
\begin{equation*}
    \mathbf{A}(\mathbf{I}+\mathbf{F}_{m})\bm{d}_{m}=\bm{r}_{m}.
\end{equation*} 
Such $\mathbf{F}_{m}$ can always be found for nonsingular $\mathbf{A}$.
Further suppose that $\|\mathbf{F}_{m}\| < 1, \; \forall m$, and hence $\mathbf{I}+\mathbf{F}_{m}$ are nonsingular.
Then,
\[ \bm{x}_{m+1} - \bm{x} = -(\mathbf{I}+\mathbf{F}_{m})^{-1}\left[(\alpha_m-1)\mathbf{I}-\mathbf{F}_{m}\right](\bm{x}_{m}-\bm{x}) .\]

Note that 
\begin{align*}
|\alpha_m-1| &= \frac{|\bm{w}_m^T \mathbf{A}\mathbf{F}_{m}\bm{d}_{m}|}{\|\bm{w}_m\|^2}
\leq     \condA
\|\mathbf{F}_{m}\| ,
\end{align*}
which implies
\begin{align*} \|\bm{x}_{m+1}-\bm{x}\| &\leq \frac{|\alpha_m-1|+\|\mathbf{F}_{m}\|}{1-\|\mathbf{F}_{m}\|}\|\bm{x}_{m}-\bm{x}\| \\
&\leq \frac{\left(1+
\condA
\right)\|\mathbf{F}_{m}\|}{1-\|\mathbf{F}_{m}\|}\|\bm{x}_{m}-\bm{x}\| .
\end{align*}
Hence, if $\|\mathbf{F}_{m}\|<(2+\condA)^{-1}$, we then have
\begin{equation*}
    (|\alpha_m-1|+\|\mathbf{F}_{m}\|) (1-\|\mathbf{F}_{m}\|)^{-1} < 1,
\end{equation*}
thus allowing us to conclude that
$\bm{x}_{m}\rightarrow \bm{x}$ and therefore our stable IR converges to the correct solution\footnote{The analysis of the error term in computing ${\bm r}_m$ is similar to \cite{Moler} and results in another dependence on $\condA$.}.

\section{Variants of stable IR} \label{sec4}

One issue with Algorithm \ref{alg:sir} arises when $\bm{r}_{m}$ is nearly orthogonal to $\mathbf{A}\bm{d}_m$, and thus $\alpha_m$ is small resulting in a minimal update. 
To address this, our stable IR algorithm can be extended to perform the line search in multiple directions. For instance, we keep track of the last $k$ directions $\bm{d}_{m}$ to form the $n \times k$ matrix $\mathbf{D}_m = [\bm{d}_{m}, \bm{d}_{m-1}, \ldots, \bm{d}_{m-k+1}]$, and keep track of the last $k$ vectors $\bm{w}_m=\mathbf{A}\bm{d}_m$ to form $\mathbf{A}\mathbf{D}_m$,
where $k \ll n$.
Solving the corresponding least squares normal equation $\bm{c}_m =\mathrm{argmin}_{\bm{c}} \|\bm{r}_{m}- \mathbf{A}\mathbf{D}_m\bm{c}\|$ (which is on the order of $O(k^3)$), we obtain the stable IR algorithm presented in Algorithm \ref{alg:sir-variant1}, where $\bm{c}_m $ is a $k$-vector and $\mathbf{D}_{m}^T\mathbf{A}^T\mathbf{A}\mathbf{D}_m$ is a relatively small $k \times k$ matrix.
\RestyleAlgo{ruled}
\begin{algorithm} 
\caption{Stable IR with multiple directions}\label{alg:sir-variant1}
\KwData{$\mathbf{A}$, $\bm{b}$, $\bm{x}_0$}
\KwResult{$\bm{x}$ such that $\mathbf{A}\bm{x}=\bm{b}$}
compute residual $\bm{r}_0 \gets \bm{b}-\mathbf{A}\bm{x}_0$\;
\While{Stopping criteria is not satisfied}{
solve for $\bm{z}$ in $\mathbf{A}\bm{z}=\bm{r}_{m}$ using basic method resulting in solution $\bm{d}_{m}$\;
construct $\mathbf{D}_m = [\bm{d}_{m}, \bm{d}_{m-1}, \ldots, \bm{d}_{m-k+1}]$\;
solve for $\bm{c}_m $ in $\mathbf{D}_m^T\mathbf{A}^T\mathbf{AD}_m\bm{c}_m = \mathbf{D}_m^T\mathbf{A}^T\bm{r}_{m}$\; 
update $\bm{x}_{m+1} \gets \bm{x}_{m} + \mathbf{D}_m \bm{c}_m$\;
update $\bm{r}_{m+1} \gets \bm{r}_{m} - \mathbf{AD}_m \bm{c}_m$\;
}
\end{algorithm}

An additional variant consists of only updating $\bm{x}_{m+1}$ in high precision after $k$ iterations of the basic method to generate a sequence of directions and compute the $k$-vector $\bm{c}_m$.
This stable IR algorithm, presented in 
Algorithm \ref{alg:sirvariant2}, is more applicable to nondeterministic basic methods (such as analog crossbar arrays) which are much faster (or power efficient) than the high-precision methods. In this case, the assumption is that the basic method is stochastic, i.e., the $\bm{d}_{m}^j$ are different vectors and $\mathbf{D}_{m}$ has rank $> 1$.
Note that if $\mathbf{A}^{-1}\bm{r}_{m}$ is in the range space of $\mathbf{D}_{m}$, then $\bm{x}_{m+1} = \mathbf{A}^{-1}\bm{b}=\bm{x}$.

\RestyleAlgo{ruled}
\begin{algorithm}
\caption{Stable IR with multiple directions II}\label{alg:sirvariant2}
\KwData{$\mathbf{A}$, $\bm{b}$, $\bm{x}_0$}
\KwResult{$\bm{x}$ such that $\mathbf{A}\bm{x}=\bm{b}$}
compute residual $\bm{r}_0 \gets \bm{b}-\mathbf{A}\bm{x}_0$\;
\While{Stopping criteria is not satisfied}{
solve for $\bm{z}$ in $\mathbf{A}\bm{z} = \bm{r}_{m}$ using the basic method $k$ times resulting in solutions $\bm{d}_{m}^j$ for $j = 1,\ldots,k$\;
construct $\mathbf{D}_{m} = [\bm{d}_{m}^1, \bm{d}_{m}^2, \ldots, \bm{d}_{m}^k]$\;
solve for $\bm{c}_m $ in $\mathbf{D}_{m}^T\mathbf{A}^T\mathbf{A}\mathbf{D}_m\bm{c}_m = \mathbf{D}_m^T\mathbf{A}^T\bm{r}_{m}$\; 
update $\bm{x}_{m+1} \gets \bm{x}_{m} + \mathbf{D}_{m} \bm{c}_m$\;
update $\bm{r}_{m+1} \gets \bm{r}_{m} - \mathbf{A}\mathbf{D}_m \bm{c}_m$\;
}
\end{algorithm}

The nondivergence and convergence results in Section \ref{sec3} also hold 
for both Algorithm \ref{alg:sir-variant1} and Algorithm \ref{alg:sirvariant2}.

\section{Experimental results}

In this section we illustrate the performance of our stable IR algorithms presented 
in Sections \ref{sec3} and \ref{sec4}. Our experiments are conducted in a MATLAB 
environment with $64$-bit arithmetic, on a single core of a computing system equipped 
with a $2.50$ GHz Octa-Core Intel Core i$7$ processor and $64$ GBs of system memory.

We implement the basic method $\mathbf{A}\bm{d}_m = \bm{r}_m$ for both standard 
IR (Algorithm \ref{alg:ir}) and stable IR (Algorithms 
\ref{alg:sir}, \ref{alg:sir-variant1}, \ref{alg:sirvariant2}) via a 
non-stationary class of iterative algorithms, known as Krylov subspace iterative solvers \cite{saad2003iterative}. Krylov subspace iterative solvers 
approximate the solution $\bm{d}_m$ from the expanding (Krylov) subspace 
$\mbox{span}\{\bm{r}_m,\mathbf{A}\bm{r}_m,\mathbf{A}^2\bm{r}_m,\ldots,\mathbf{A}^j\bm{r}_m\}$, 
where $j$ denotes the iteration index. For the sake of generality, in this paper 
we only consider general-purpose solvers such as GMRES, FGMRES, CGS, and BiCGSTAB, but an extensive list can be found in  
\cite{liesen2013krylov}. For each one of these methods, the most expensive operation is the computation of matrix-vector products with the matrix $\mathbf{A}$. Therefore, reducing the wall-clock time of Krylov subspace iterative solvers requires the faster computation of matrix-vector products with the matrix $\mathbf{A}$. 
Throughout the rest 
of our numerical experiments section we assume that these matrix-vector products 
are performed via: $a$) a simulated analog crossbar array; and $b$) the use of 
reduced precision digital arithmetic. 

\subsection{Analog in-memory crossbar array} \label{sec:analog}

We begin with the consideration of a model analog hardware accelerator equipped with a crossbar array of resistive processing units (RPUs). To perform the simulation of the analog hardware, we used a MATLAB version of the publicly available simulator~\cite{rasch2021flexible} with a PyTorch interface for emulating the noise, timing, energy and ADC/DAC characteristics of an analog crossbar array. This simulator models all sources of analog noise as scaled Gaussian processes and are based on currently realizable analog hardware~\cite{Gokmen2016}, where we used the same parameters as in \cite{kalantzis2021solving}.  

\subsubsection{Random matrices}
We first consider a symmetric and positive-definite matrix $\mathbf{A}$ of size $\Mu=2000$, 
whose entries are set as 
\begin{equation*} 
    [\mathbf{A}]_{ij} =
    \begin{cases}
      1+\sqrt{i}                     & \mathtt{if}\ \ \ i==j \\
      1/|i-j|        & \mathtt{if}\ \ \ i\neq j
    \end{cases}.
\end{equation*}
The decay in the off-diagonal entries of $\bA$ models a progressively 
decaying correlation among the feature space of a signal or general data collection. 
Model covariance matrices of this form have already been considered in the context of 
IR \cite{LeGallo,kalantzis2013accelerating,kalantzis2018scalable}.

\begin{figure}[htbp]
\centerline{\includegraphics[width=0.23\textwidth]{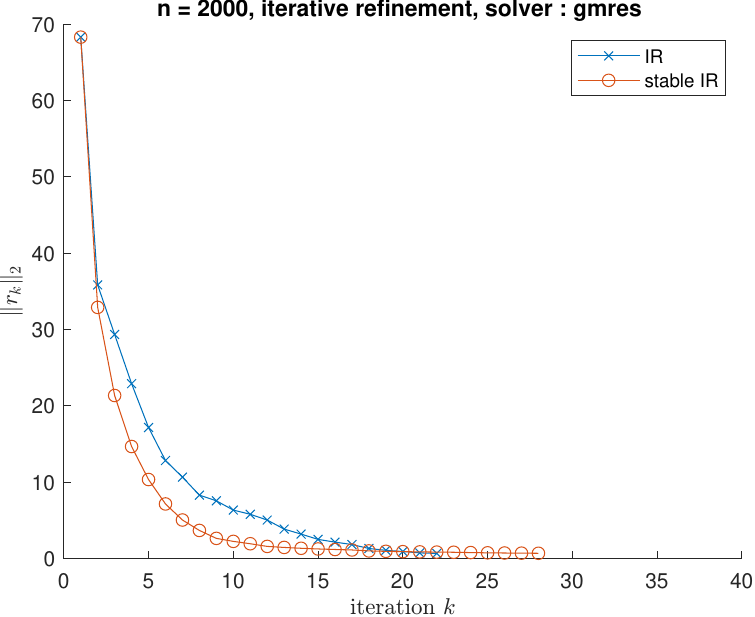}
\quad\includegraphics[width=0.23\textwidth]{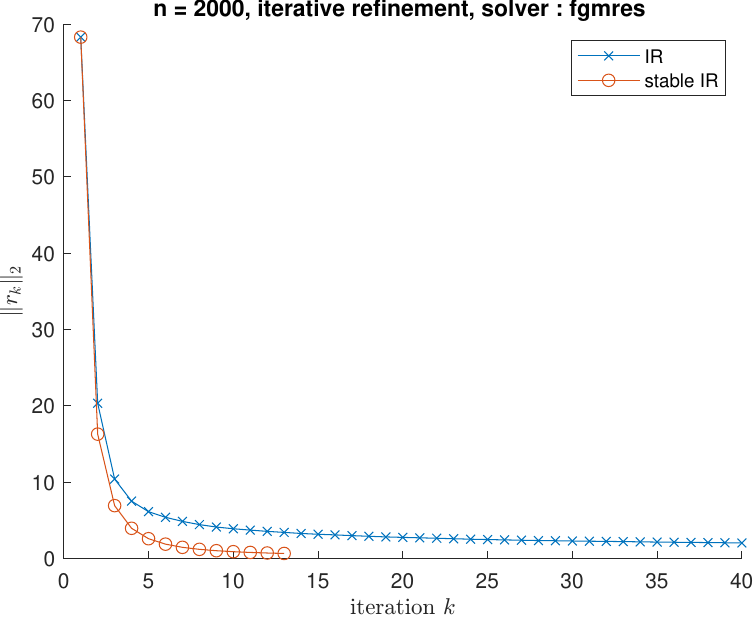}}
\centerline{(a)\hspace{0.23\textwidth}(b)}
\caption{Solving $\mathbf{A}\bm{x}=\bm{b}$ using RPU. $\mathbf{A}$ is a symmetric matrix with decreasing off-diagonal elements. (a) GMRES-IR. (b) FGMRES-IR. }\label{fig:rpu-gmres}
\end{figure}

\begin{figure}[htbp]
\centerline{\includegraphics[width=0.23\textwidth]{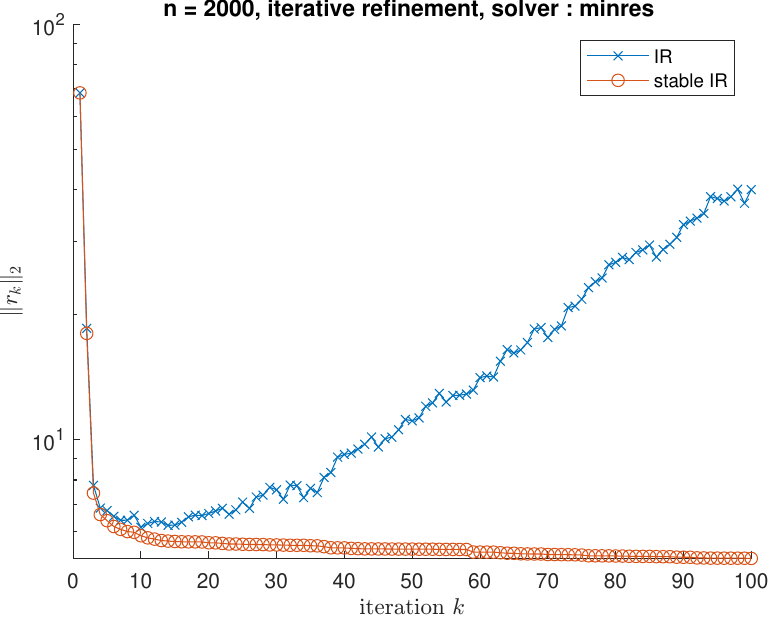}
\quad\includegraphics[width=0.23\textwidth]{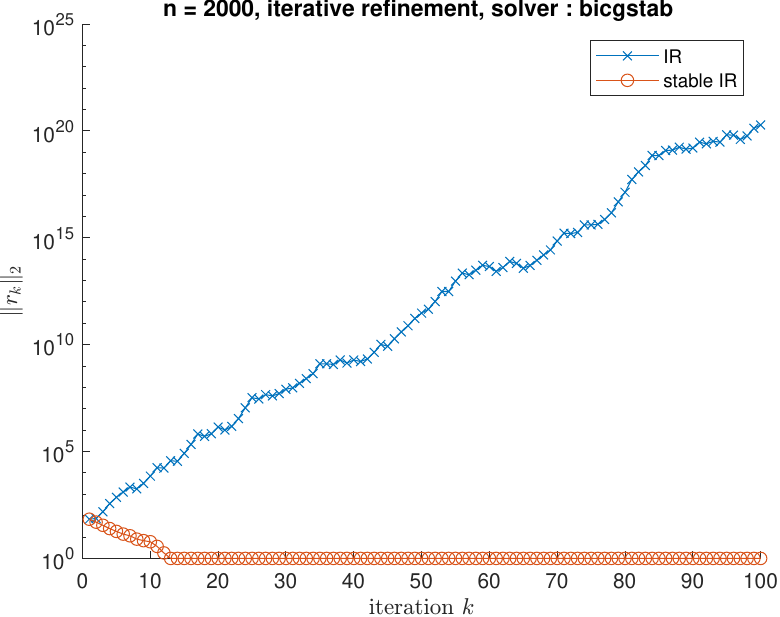}}
\centerline{(a)\hspace{0.2\textwidth}(b)}
\caption{Solving $\mathbf{A}\bm{x}=\bm{b}$ using RPU. $\mathbf{A}$ is symmetric with decreasing off-diagonal elements. (a) MINRES-IR. (b) BiCGSTAB-IR. }\label{fig:rpu-minres}
\end{figure}

\begin{figure}[htbp]
\centerline{\includegraphics[width=0.23\textwidth]{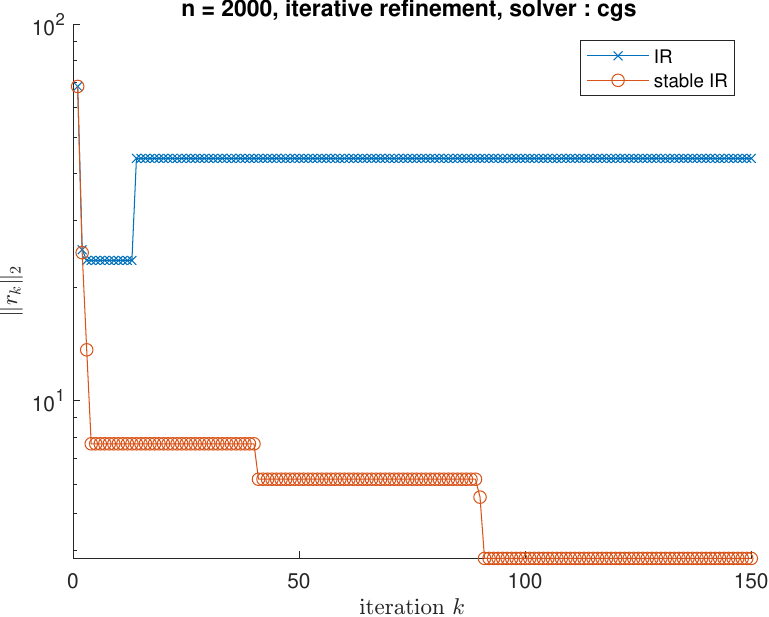}
\quad\includegraphics[width=0.23\textwidth]{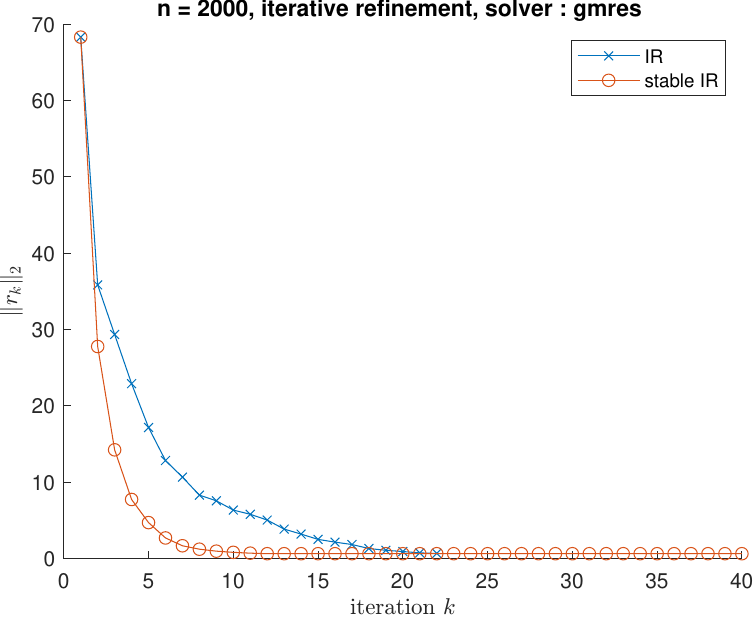}}
\centerline{(a)\hspace{0.2\textwidth}(b)}
\caption{Solving $\mathbf{A}\bm{x}=\bm{b}$ using RPU. $\mathbf{A}$ is symmetric with decreasing off-diagonal elements. (a) CGS-IR.  (b) GMRES-IR running the basic method $k = 10$ times between each update (Algorithm \ref{alg:sirvariant2}).}\label{fig:rpu-cgs}
\end{figure}

\begin{figure}[htbp]
\centerline{\includegraphics[width=0.23\textwidth]{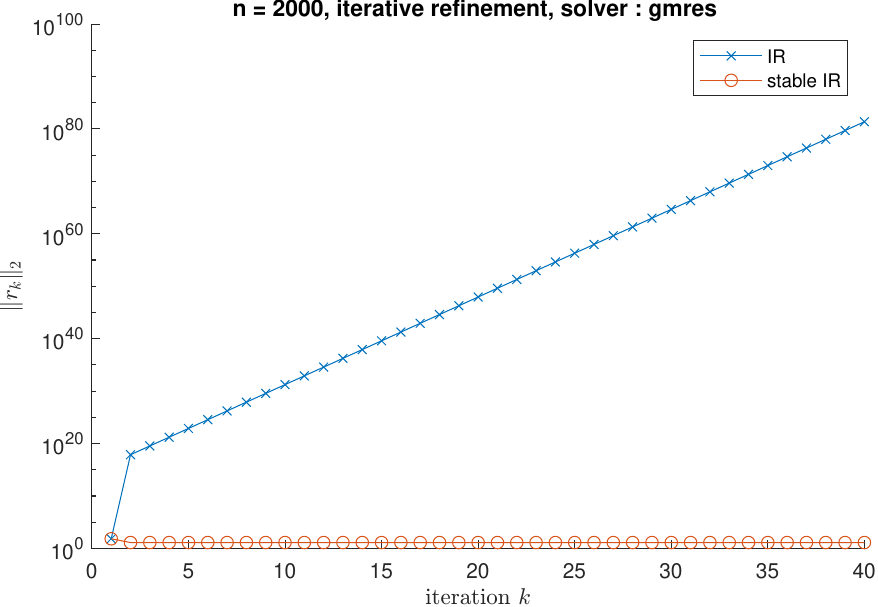}
\quad\includegraphics[width=0.23\textwidth]{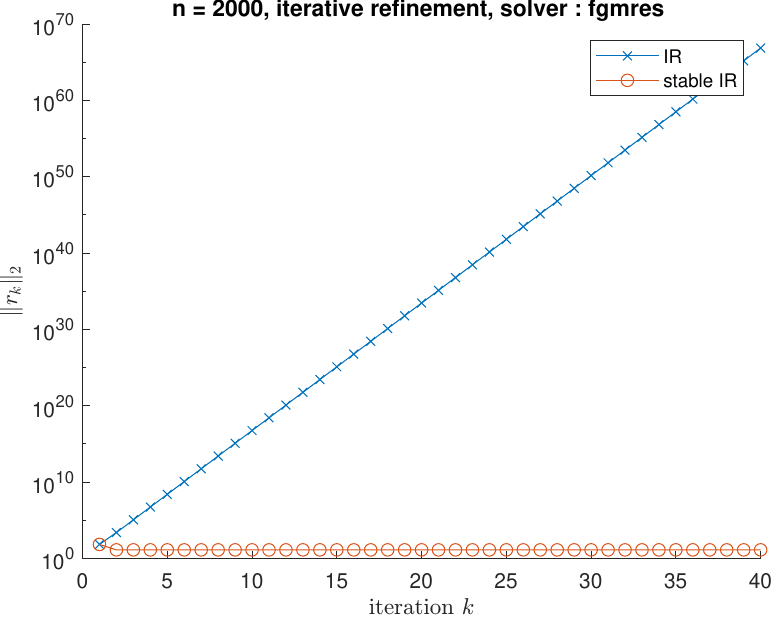}}
\centerline{(a)\hspace{0.23\textwidth}(b)}
\caption{Solving $\mathbf{A}\bm{x}=\bm{b}$ using RPU. $\mathbf{A}$ is a random $2000 \times 2000$ matrix of i.i.d.\ entries from a uniform distribution. (a) GMRES-IR. (b) FGMRES-IR.}\label{fig:rpu-gmres-random}
\end{figure}

Figs.~\ref{fig:rpu-gmres},~\ref{fig:rpu-minres} and~\ref{fig:rpu-cgs}(a) plot the residual norm ($y$-axis) of classical and stable IR as the number of iterations increases ($x$-axis). We see 
that the performance of our stable IR algorithm is superior to that of the classical IR 
algorithm, with the classical IR algorithm diverging when the basic method is MINRES 
or BiCGSTAB, whereas the residual norm remains small under our stable IR algorithm.
Next, we consider a different matrix $\mathbf{A}$ of size $\Mu=2000$ 
where each entry $[\mathbf{A}]_{ij}$ is taken i.i.d.\ from a uniform 
distribution in $[0,1]$. In this case, $\mathbf{A}$ is generally 
neither positive definite nor symmetric. In Fig.~\ref{fig:rpu-gmres-random}, we observe that GMRES-IR and FGMRES-IR diverge, while the stable GMRES-IR and stable FGMRES-IR algorithms do not.

\subsubsection{Matrices from real-world examples}

\begin{figure}[htbp]
\centerline{\includegraphics[width=0.23\textwidth]{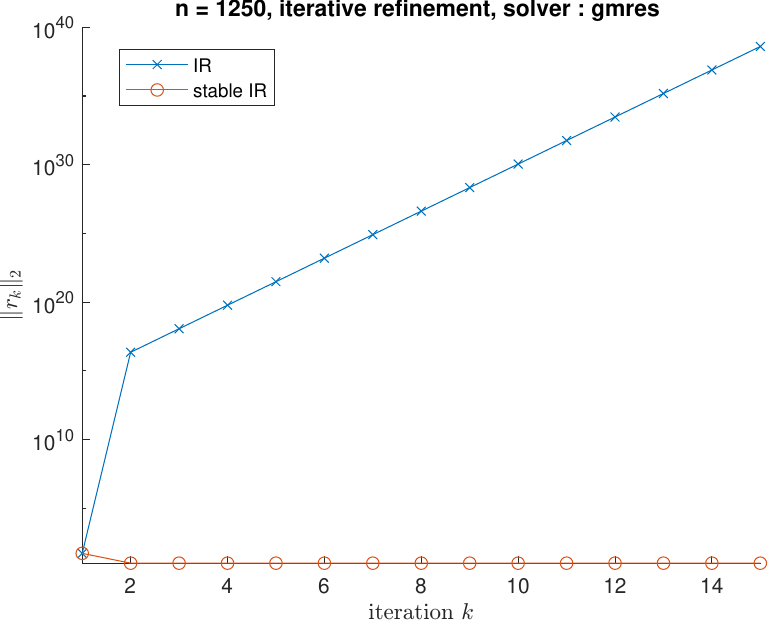}
\quad\includegraphics[width=0.23\textwidth]{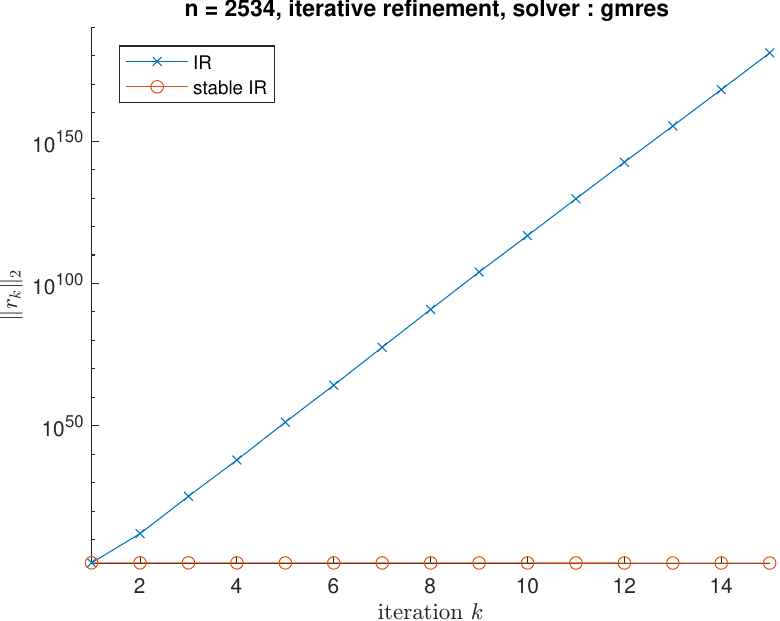}}
\centerline{(a)\hspace{0.23\textwidth}(b)}
\caption{Solving $\mathbf{A}\bm{x}=\bm{b}$ using RPU and GMRES-IR. (a) $\mathbf{A}$ is from a fluid dynamics problem. (b) $\mathbf{A}$ is from an electromagnetics problem.}
\label{fig:rpu-gmres-sparse}
\end{figure}

Next, we consider several matrices from real-world examples taken from the SuiteSparse Matrix Collection \cite{suitesparse2011}, normalized to the dynamic range of the analog array. In particular, in Fig.~\ref{fig:rpu-gmres-sparse} we show (a) a matrix of size $\Mu=1250$ from a reaction-diffusion Brusselator model (\texttt{rdb1250l}) and (b) a matrix of size $\Mu=2534$ from an electromagnetic field model (\texttt{qc2534}).

\begin{figure}[htbp]
\centerline{\includegraphics[width=0.23\textwidth]{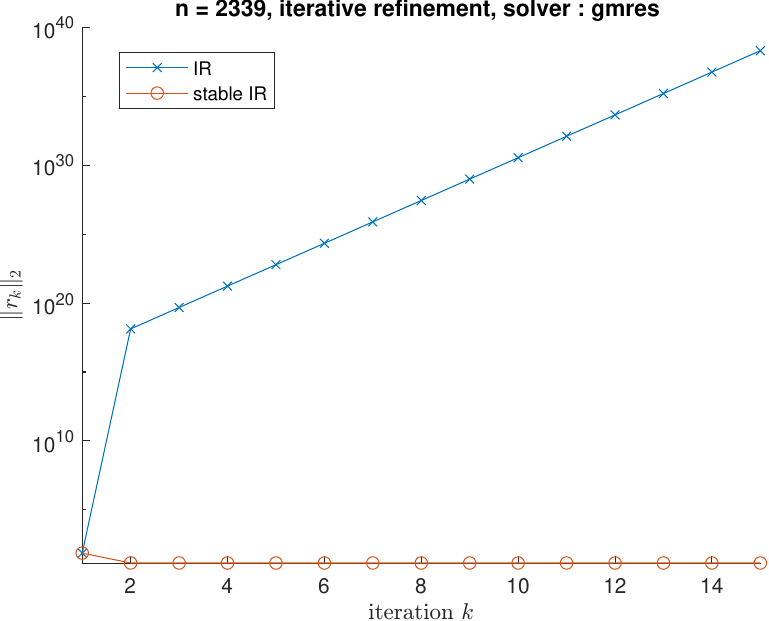}
\quad\includegraphics[width=0.23\textwidth]{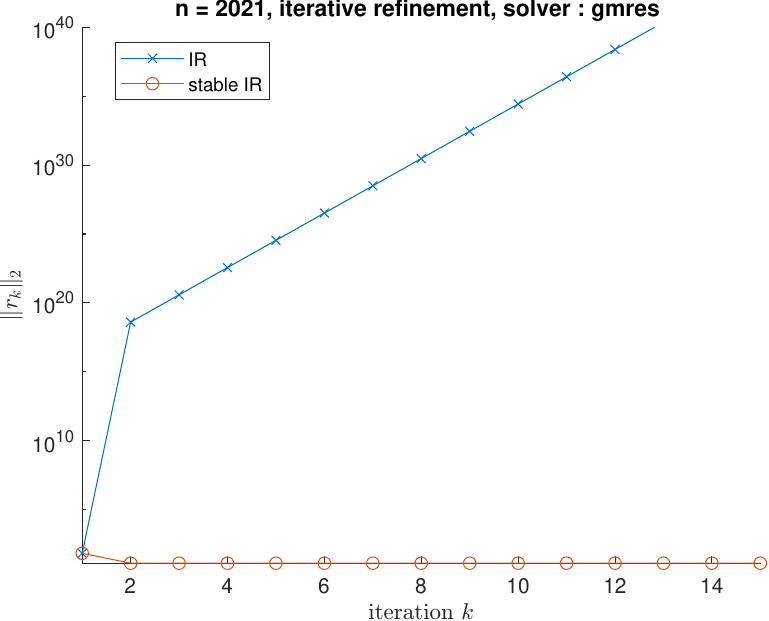}}
\centerline{(a)\hspace{0.23\textwidth}(b)}
\caption{Solving $\mathbf{A}\bm{x}=\bm{b}$ using RPU and GMRES-IR. (a) $\mathbf{A}$ is from a finite element analysis of the heart. (b) $\mathbf{A}$ is from a chemical process simulation problem.}
\label{fig:rpu-gmres-sparse2}
\end{figure}

In Fig.~\ref{fig:rpu-gmres-sparse2} we consider (a) a matrix of size $\Mu=2339$ from a finite element model of a heart (\texttt{heart3}) and (b) a matrix of size $\Mu=2021$ from a chemical process simulation model  (\texttt{west2021}).
We see that in all these cases, the standard IR algorithm diverges while the residual error in the proposed stable IR algorithm remain small.

\subsection{Low-precision digital accelerator}

\begin{figure}[htbp]
\centerline{\includegraphics[width=0.23\textwidth]{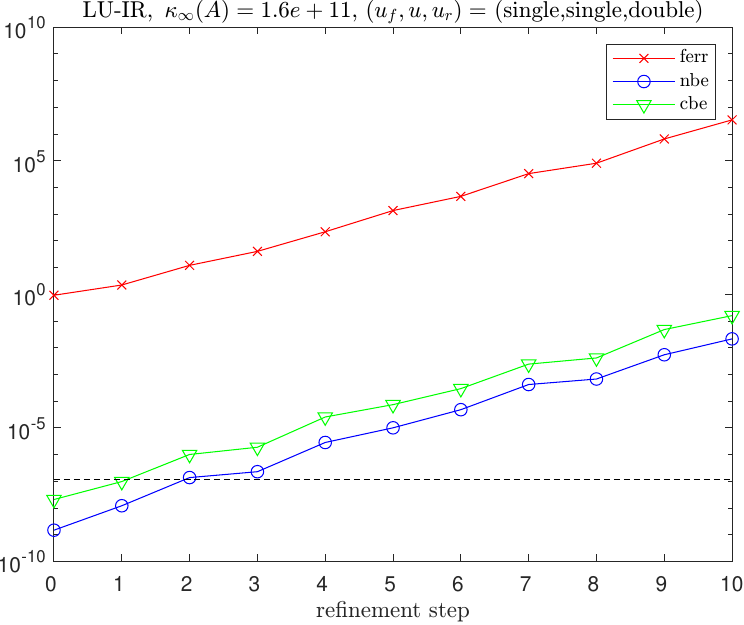}
\quad\includegraphics[width=0.23\textwidth]{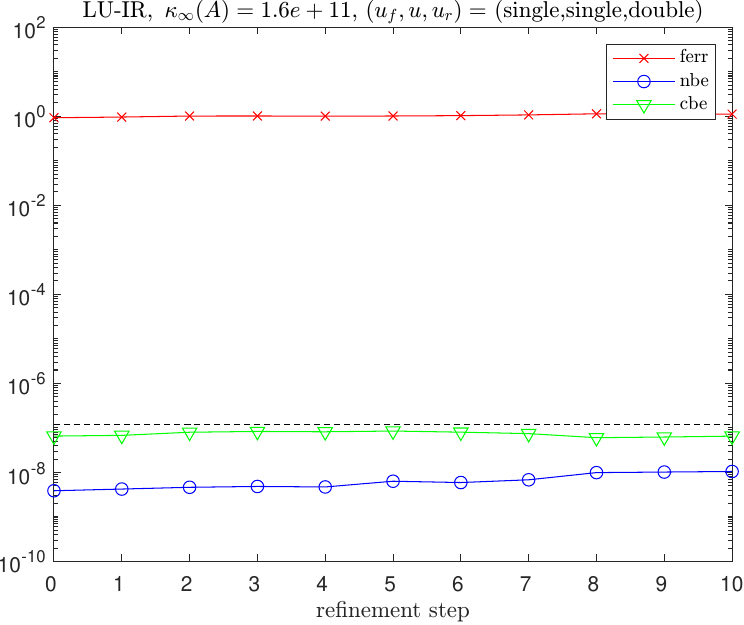}}
\centerline{(a)\hspace{0.23\textwidth}(b)}
\caption{LU-IR. \texttt{ferr}, \texttt{nbe} and \texttt{cbe} are the forward, normwise backward and componentwise backward error as defined in \cite{IR3precision:2018}. (a) IR (Algorithm \ref{alg:ir}). (b) Stable IR (Algorithm \ref{alg:sir}).}\label{fig:LU-IR}
\end{figure}

Our stable IR algorithms will also benefit from scenarios where the matrix-vector products during the iterative solution of the basic method are performed in lower 
digital precision while the rest of the steps are still performed in higher digital precision.
For instance, consider the example in~\cite{IR3precision:2018} where the basic method 
for solving $\mathbf{A}\bm{x}=\bm{b}$ is LU decomposition in single precision, single precision is used to store intermediate results, and the residual and $\alpha_m$ (if applicable) are computed in double precision. 
For this scheme, upon adapting the code provided in \cite{IR3precision:2018},  Fig.~\ref{fig:LU-IR}(a) shows that
LU-IR diverges for both $\|\bm{x}_{m}-\bm{x}\|$ and $\|\bm{r}_{m}\|$ under a matrix $\mathbf{A}$
with $\condA = 1.6\times 10^{11}$.
On the other hand, Fig.~\ref{fig:LU-IR}(b) shows that for the stable LU-IR scheme (Algorithm \ref{alg:sir}), $\bm{x}_{m}$ are bounded and the residual errors are within the desired precision.

\subsection{Variants of stable IR}

\begin{figure}[htbp]
\centerline{\includegraphics[width=0.23\textwidth]{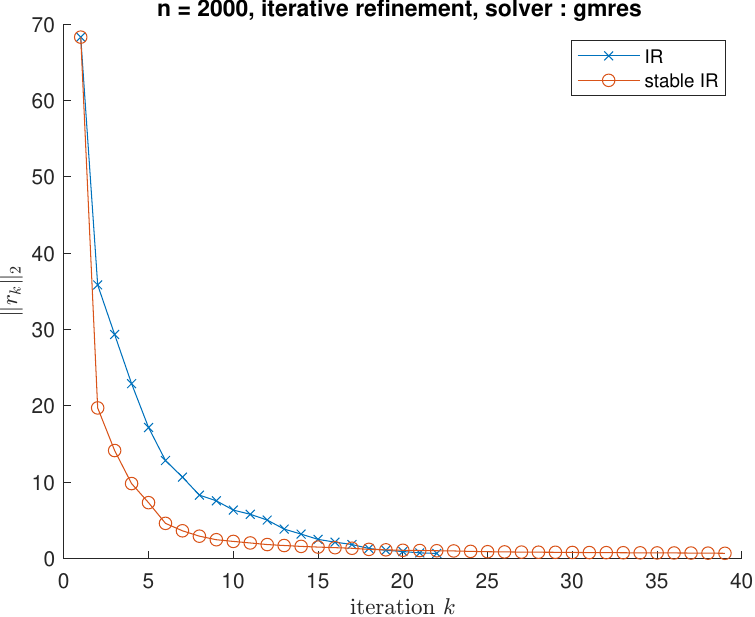}
\quad\includegraphics[width=0.23\textwidth]{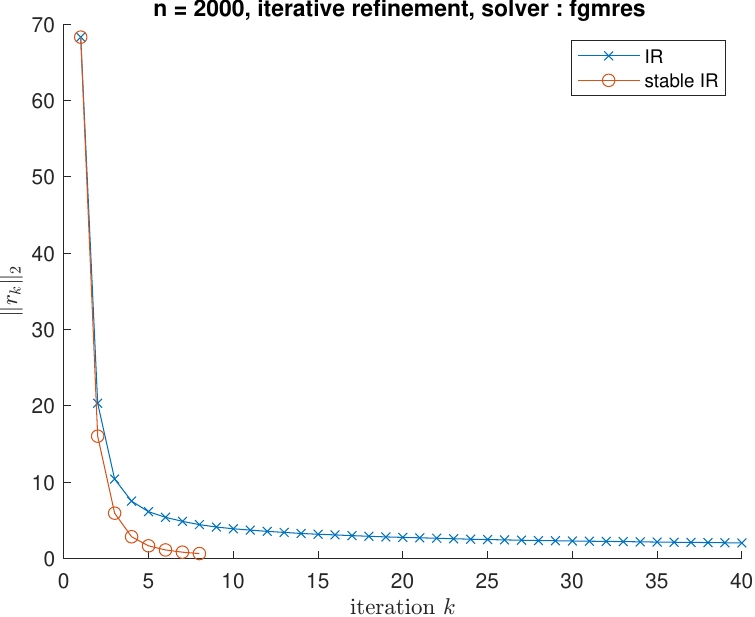}}
\centerline{(a)\hspace{0.23\textwidth}(b)}
\caption{Solving $\mathbf{A}\bm{x}=\bm{b}$ using RPU with $k = 10$ previous directions. $\mathbf{A}$ is a symmetric matrix with decreasing off-diagonal elements. (a) GMRES. (b) FGMRES. }\label{fig:rpu-gmres-k-10}
\end{figure}

Applying Algorithm \ref{alg:sir-variant1} with $k=10$ to the setup in Section \ref{sec:analog} 
shows in Fig.~\ref{fig:rpu-gmres-k-10} that it performs better than Algorithm \ref{alg:sir} with $k=1$ (Fig.~\ref{fig:rpu-gmres}) for both GMRES-IR and FGMRES-IR, converging to the solution in much less iterations than Algorithm \ref{alg:ir}. 
Similarly, for Algorithm~\ref{alg:sirvariant2} with $k=10$, we observe in 
Fig.~\ref{fig:rpu-cgs}(b)
that Algorithm~\ref{alg:sirvariant2} for GMRES-IR performs better than both Algorithm \ref{alg:sir} in Fig.~\ref{fig:rpu-gmres}(a) and Algorithm~\ref{alg:sir-variant1} in Fig.~\ref{fig:rpu-gmres-k-10}(a).

\section{Conclusions}

In this paper we presented a modification of the classical
IR
algorithm that guarantees nondivergence under unknown or arbitrarily large error for the basic method and that demonstrates performance superior to IR in practice. This is especially useful for emerging inaccurate computing paradigms where the variance of the noise/error can be unknown and/or large. 
One issue of the stable IR algorithm is that when $\bm{r}_{m}$ is nearly orthogonal to $\mathbf{A}\bm{d}_m$, this results in $\alpha_m$ being small resulting in minimal update. As part of our future work, we aim to find an alternative direction for $\bm{d}_{m}$ with which to update $\bm{x}_{m}$. 

\vfill\pagebreak



\begin{thebibliography}{10}
\providecommand{\url}[1]{#1}
\csname url@samestyle\endcsname
\providecommand{\newblock}{\relax}
\providecommand{\bibinfo}[2]{#2}
\providecommand{\BIBentrySTDinterwordspacing}{\spaceskip=0pt\relax}
\providecommand{\BIBentryALTinterwordstretchfactor}{4}
\providecommand{\BIBentryALTinterwordspacing}{\spaceskip=\fontdimen2\font plus
\BIBentryALTinterwordstretchfactor\fontdimen3\font minus
  \fontdimen4\font\relax}
\providecommand{\BIBforeignlanguage}[2]{{%
\expandafter\ifx\csname l@#1\endcsname\relax
\typeout{** WARNING: IEEEtran.bst: No hyphenation pattern has been}%
\typeout{** loaded for the language `#1'. Using the pattern for}%
\typeout{** the default language instead.}%
\else
\language=\csname l@#1\endcsname
\fi
#2}}
\providecommand{\BIBdecl}{\relax}
\BIBdecl

\bibitem{markovsky2007overview}
I.~Markovsky and S.~Van~Huffel, ``Overview of total least-squares methods,''
  \emph{Signal processing}, vol.~87, no.~10, pp. 2283--2302, 2007.

\bibitem{cadzow1990signal}
J.~A. Cadzow, ``Signal processing via least squares error modeling,''
  \emph{IEEE ASSP Magazine}, vol.~7, no.~4, pp. 12--31, 1990.

\bibitem{selesnick2013least}
I.~Selesnick, ``Least squares with examples in signal processing,''
  \emph{Connexions}, vol.~4, pp. 1--25, 2013.

\bibitem{wang2011least}
D.~Wang and F.~Ding, ``Least squares based and gradient based iterative
  identification for {W}iener nonlinear systems,'' \emph{Signal Processing},
  vol.~91, no.~5, pp. 1182--1189, 2011.

\bibitem{aastrom1971system}
K.~J. {\AA}str{\"o}m and P.~Eykhoff, ``System identification—a survey,''
  \emph{Automatica}, vol.~7, no.~2, pp. 123--162, 1971.

\bibitem{lu2014overview}
L.~Lu \emph{et~al.}, ``An overview of massive {MIMO}: Benefits and
  challenges,'' \emph{IEEE journal of selected topics in signal processing},
  vol.~8, no.~5, pp. 742--758, 2014.

\bibitem{marzetta2015massive}
T.~L. Marzetta, ``Massive {MIMO}: an introduction,'' \emph{Bell Labs Technical
  Journal}, vol.~20, pp. 11--22, 2015.

\bibitem{larsson2014massive}
E.~G. Larsson, O.~Edfors, F.~Tufvesson, and T.~L. Marzetta, ``Massive {MIMO}
  for next generation wireless systems,'' \emph{IEEE communications magazine},
  vol.~52, no.~2, pp. 186--195, 2014.

\bibitem{saad2003iterative}
Y.~Saad, \emph{Iterative methods for sparse linear systems}.\hskip 1em plus
  0.5em minus 0.4em\relax SIAM, 2003.

\bibitem{zhao2000hierarchical}
M.~Zhao, R.~V. Panda, S.~S. Sapatnekar, T.~Edwards, R.~Chaudhry, and D.~Blaauw,
  ``Hierarchical analysis of power distribution networks,'' in
  \emph{Proceedings of the 37th Annual Design Automation Conference}, 2000, pp.
  150--155.

\bibitem{chang2000analysis}
P.~S. Chang and A.~N. Willson, ``Analysis of conjugate gradient algorithms for
  adaptive filtering,'' \emph{IEEE Transactions on Signal Processing}, vol.~48,
  no.~2, pp. 409--418, 2000.

\bibitem{tu2020efficient}
J.~Tu, M.~Lou, J.~Jiang, D.~Shu, and G.~He, ``An efficient massive {MIMO}
  detector based on second-order richardson iteration: From algorithm to
  flexible architecture,'' \emph{IEEE Transactions on Circuits and Systems I:
  Regular Papers}, vol.~67, no.~11, pp. 4015--4028, 2020.

\bibitem{wang2022randomized}
Z.~Wang, R.~M. Gower, Y.~Xia, L.~He, and Y.~Huang, ``Randomized iterative
  methods for low-complexity large-scale {MIMO} detection,'' \emph{IEEE
  Transactions on Signal Processing}, vol.~70, pp. 2934--2949, 2022.

\bibitem{liesen2013krylov}
J.~Liesen and Z.~Strakos, \emph{Krylov subspace methods: principles and
  analysis}.\hskip 1em plus 0.5em minus 0.4em\relax Oxford University Press,
  2013.

\bibitem{Wilkinson}
J.~H. Wilkinson, \emph{Rounding errors in algebraic processes}.\hskip 1em plus
  0.5em minus 0.4em\relax SIAM, 2023.

\bibitem{golub1999inexact}
G.~H. Golub and Q.~Ye, ``Inexact preconditioned conjugate gradient method with
  inner-outer iteration,'' \emph{SIAM Journal on Scientific Computing},
  vol.~21, no.~4, pp. 1305--1320, 1999.

\bibitem{gleich2010inner}
D.~F. Gleich, A.~P. Gray, C.~Greif, and T.~Lau, ``An inner-outer iteration for
  computing pagerank,'' \emph{SIAM Journal on Scientific Computing}, vol.~32,
  no.~1, pp. 349--371, 2010.

\bibitem{haidar2017investigating}
A.~Haidar, P.~Wu, S.~Tomov, and J.~Dongarra, ``Investigating half precision
  arithmetic to accelerate dense linear system solvers,'' in \emph{Proceedings
  of the 8th workshop on latest advances in scalable algorithms for large-scale
  systems}, 2017, pp. 1--8.

\bibitem{Moler}
C.~B. Moler, ``Iterative refinement in floating point,'' \emph{Journal of the
  ACM (JACM)}, vol.~14, no.~2, pp. 316--321, 1967.

\bibitem{haidar2018harnessing}
A.~Haidar, S.~Tomov, J.~Dongarra, and N.~J. Higham, ``Harnessing {GPU} tensor
  cores for fast {FP16} arithmetic to speed up mixed-precision iterative
  refinement solvers,'' in \emph{SC18: International Conference for High
  Performance Computing, Networking, Storage and Analysis}.\hskip 1em plus
  0.5em minus 0.4em\relax IEEE, 2018, pp. 603--613.

\bibitem{haidar2020mixed}
A.~Haidar, H.~Bayraktar, S.~Tomov, J.~Dongarra, and N.~J. Higham,
  ``Mixed-precision iterative refinement using tensor cores on {GPU}s to
  accelerate solution of linear systems,'' \emph{Proceedings of the Royal
  Society A}, vol. 476, no. 2243, p. 20200110, 2020.

\bibitem{haidar2018design}
A.~Haidar \emph{et~al.}, ``The design of fast and energy-efficient linear
  solvers: On the potential of half-precision arithmetic and iterative
  refinement techniques,'' in \emph{International conference on computational
  science}.\hskip 1em plus 0.5em minus 0.4em\relax Springer, 2018, pp.
  586--600.

\bibitem{kalantzis2021solving}
V.~Kalantzis \emph{et~al.}, ``Solving sparse linear systems with approximate
  inverse preconditioners on analog devices,'' in \emph{2021 IEEE High
  Performance Extreme Computing Conference (HPEC)}.\hskip 1em plus 0.5em minus
  0.4em\relax IEEE, 2021, pp. 1--7.

\bibitem{kalantzis2023fgmres}
------, ``Solving sparse linear systems via flexible {GMRES} with in-memory
  analog preconditioning,'' in \emph{2023 IEEE High Performance Extreme
  Computing Conference (HPEC)}.\hskip 1em plus 0.5em minus 0.4em\relax IEEE,
  2023.

\bibitem{stochastic_computing_survey}
\BIBentryALTinterwordspacing
A.~Alaghi and J.~P. Hayes, ``Survey of stochastic computing,'' \emph{ACM Trans.
  Embed. Comput. Syst.}, vol.~12, no.~2s, may 2013. [Online]. Available:
  \url{https://doi.org/10.1145/2465787.2465794}
\BIBentrySTDinterwordspacing

\bibitem{dither_computing}
C.~Wu, ``Dither computing: a hybrid deterministic-stochastic computing
  framework,'' in \emph{2021 IEEE 28th Symposium on Computer Arithmetic
  (ARITH)}, 2021, pp. 70--77.

\bibitem{rasch2021flexible}
M.~J. Rasch, D.~Moreda, T.~Gokmen, M.~L. Gallo, F.~Carta, C.~Goldberg, K.~E.
  Maghraoui, A.~Sebastian, and V.~Narayanan, ``A flexible and fast {PyTorch}
  toolkit for simulating training and inference on analog crossbar arrays,''
  \emph{arXiv preprint arXiv:2104.02184}, 2021.

\bibitem{Gokmen2016}
T.~Gokmen and Y.~Vlasov, ``Acceleration of deep neural network training with
  resistive cross-point devices: Design considerations,'' \emph{Frontiers in
  Neuroscience}, vol.~10, 2016.

\bibitem{LeGallo}
M.~Le~Gallo \emph{et~al.}, ``Mixed-precision in-memory computing,''
  \emph{Nature Electronics}, vol.~1, no.~4, pp. 246--253, 2018.

\bibitem{kalantzis2013accelerating}
``Accelerating data uncertainty quantification by solving linear systems with
  multiple right-hand sides,'' \emph{Numerical Algorithms}, vol.~62, pp.
  637--653, 2013.

\bibitem{kalantzis2018scalable}
V.~Kalantzis, A.~C.~I. Malossi, C.~Bekas, A.~Curioni, E.~Gallopoulos, and
  Y.~Saad, ``A scalable iterative dense linear system solver for multiple
  right-hand sides in data analytics,'' \emph{Parallel Computing}, vol.~74, pp.
  136--153, 2018.

\bibitem{suitesparse2011}
T.~A. Davis and Y.~Hu, ``The {U}niversity of {F}lorida sparse matrix
  collection,'' \emph{ACM Transactions on Mathematical Software}, vol.~38,
  no.~1, pp. 1--25, Nov. 2011.

\bibitem{IR3precision:2018}
E.~Carson and N.~J. Higham, ``Accelerating the solution of linear systems by
  iterative refinement in three precisions,'' \emph{SIAM Journal on Scientific
  Computing}, vol.~40, no.~2, pp. A817--A847, 2018.

\end{thebibliography}
\end{document}